\input amstex
\input amsppt.sty

\magnification=\magstep1 

\vcorrection{-0.8cm}

\catcode`\@=11
\font\tenln    = line10
\font\tenlnw   = linew10

\newskip\Einheit \Einheit=0.5cm
\newcount\xcoord \newcount\ycoord
\newdimen\xdim \newdimen\ydim \newdimen\PfadD@cke \newdimen\Pfadd@cke

%%%%%%%%%%%%%%%%%%%%%%%%%%%%%%%%%%%%%%%%%%%%%%%%%
%LaTeX counters, dimensions, variables for lines%
%%%%%%%%%%%%%%%%%%%%%%%%%%%%%%%%%%%%%%%%%%%%%%%%%
\newcount\@tempcnta
\newcount\@tempcntb

\newdimen\@tempdima
\newdimen\@tempdimb

\newdimen\@wholewidth
\newdimen\@halfwidth

\newcount\@xarg
\newcount\@yarg
\newcount\@yyarg
\newbox\@linechar
\newbox\@tempboxa
\newdimen\@linelen
\newdimen\@clnwd
\newdimen\@clnht

\newif\if@negarg

\def\@whilenoop#1{}
\def\@whiledim#1\do #2{\ifdim #1\relax#2\@iwhiledim{#1\relax#2}\fi}
\def\@iwhiledim#1{\ifdim #1\let\@nextwhile=\@iwhiledim
        \else\let\@nextwhile=\@whilenoop\fi\@nextwhile{#1}}

\def\@whileswnoop#1\fi{}
\def\@whilesw#1\fi#2{#1#2\@iwhilesw{#1#2}\fi\fi}
\def\@iwhilesw#1\fi{#1\let\@nextwhile=\@iwhilesw
         \else\let\@nextwhile=\@whileswnoop\fi\@nextwhile{#1}\fi}

\def\thinlines{\let\@linefnt\tenln \let\@circlefnt\tencirc
  \@wholewidth\fontdimen8\tenln \@halfwidth .5\@wholewidth}
\def\thicklines{\let\@linefnt\tenlnw \let\@circlefnt\tencircw
  \@wholewidth\fontdimen8\tenlnw \@halfwidth .5\@wholewidth}
\thinlines
%%%%%%%%%%%%%%%%%%%%%%%%%%%%%%%%%%%%%%%%%%%%%%%%%%%%%%%%%%%

\PfadD@cke1pt \Pfadd@cke0.5pt
\def\PfadDicke#1{\PfadD@cke#1 \divide\PfadD@cke by2 \Pfadd@cke\PfadD@cke \multiply\PfadD@cke by2}
\long\def\LOOP#1\REPEAT{\def\BODY{#1}\ITERATE}
\def\ITERATE{\BODY \let\next\ITERATE \else\let\next\relax\fi \next}
\let\REPEAT=\fi
\def\Punkt{\hbox{\raise-2pt\hbox to0pt{\hss$\ssize\bullet$\hss}}}
\def\DuennPunkt(#1,#2){\unskip
  \raise#2 \Einheit\hbox to0pt{\hskip#1 \Einheit
          \raise-2.5pt\hbox to0pt{\hss$\bullet$\hss}\hss}}
\def\NormalPunkt(#1,#2){\unskip
  \raise#2 \Einheit\hbox to0pt{\hskip#1 \Einheit
          \raise-3pt\hbox to0pt{\hss\twelvepoint$\bullet$\hss}\hss}}
\def\DickPunkt(#1,#2){\unskip
  \raise#2 \Einheit\hbox to0pt{\hskip#1 \Einheit
          \raise-4pt\hbox to0pt{\hss\fourteenpoint$\bullet$\hss}\hss}}
\def\Kreis(#1,#2){\unskip
  \raise#2 \Einheit\hbox to0pt{\hskip#1 \Einheit
          \raise-4pt\hbox to0pt{\hss\fourteenpoint$\circ$\hss}\hss}}

%%%%%%%%%%%%%%%%%%%%%
%LaTeX line macros%
%%%%%%%%%%%%%%%%%%%%%
\def\Line@(#1,#2)#3{\@xarg #1\relax \@yarg #2\relax
\@linelen=#3\Einheit
\ifnum\@xarg =0 \@vline
  \else \ifnum\@yarg =0 \@hline \else \@sline\fi
\fi}

\def\@sline{\ifnum\@xarg< 0 \@negargtrue \@xarg -\@xarg \@yyarg -\@yarg
  \else \@negargfalse \@yyarg \@yarg \fi
\ifnum \@yyarg >0 \@tempcnta\@yyarg \else \@tempcnta -\@yyarg \fi
\ifnum\@tempcnta>6 \@badlinearg\@tempcnta0 \fi
\ifnum\@xarg>6 \@badlinearg\@xarg 1 \fi
\setbox\@linechar\hbox{\@linefnt\@getlinechar(\@xarg,\@yyarg)}%
\ifnum \@yarg >0 \let\@upordown\raise \@clnht\z@
   \else\let\@upordown\lower \@clnht \ht\@linechar\fi
\@clnwd=\wd\@linechar
\if@negarg \hskip -\wd\@linechar \def\@tempa{\hskip -2\wd\@linechar}\else
     \let\@tempa\relax \fi
\@whiledim \@clnwd <\@linelen \do
  {\@upordown\@clnht\copy\@linechar
   \@tempa
   \advance\@clnht \ht\@linechar
   \advance\@clnwd \wd\@linechar}%
\advance\@clnht -\ht\@linechar
\advance\@clnwd -\wd\@linechar
\@tempdima\@linelen\advance\@tempdima -\@clnwd
\@tempdimb\@tempdima\advance\@tempdimb -\wd\@linechar
\if@negarg \hskip -\@tempdimb \else \hskip \@tempdimb \fi
\multiply\@tempdima \@m
\@tempcnta \@tempdima \@tempdima \wd\@linechar \divide\@tempcnta \@tempdima
\@tempdima \ht\@linechar \multiply\@tempdima \@tempcnta
\divide\@tempdima \@m
\advance\@clnht \@tempdima
\ifdim \@linelen <\wd\@linechar
   \hskip \wd\@linechar
  \else\@upordown\@clnht\copy\@linechar\fi}

\def\@hline{\ifnum \@xarg <0 \hskip -\@linelen \fi
\vrule height\Pfadd@cke width \@linelen depth\Pfadd@cke
\ifnum \@xarg <0 \hskip -\@linelen \fi}

\def\@getlinechar(#1,#2){\@tempcnta#1\relax\multiply\@tempcnta 8
\advance\@tempcnta -9 \ifnum #2>0 \advance\@tempcnta #2\relax\else
\advance\@tempcnta -#2\relax\advance\@tempcnta 64 \fi
\char\@tempcnta}

\def\Vektor(#1,#2)#3(#4,#5){\unskip\leavevmode
  \xcoord#4\relax \ycoord#5\relax
      \raise\ycoord \Einheit\hbox to0pt{\hskip\xcoord \Einheit
         \Vector@(#1,#2){#3}\hss}}

\def\Vector@(#1,#2)#3{\@xarg #1\relax \@yarg #2\relax
\@tempcnta \ifnum\@xarg<0 -\@xarg\else\@xarg\fi
\ifnum\@tempcnta<5\relax
\@linelen=#3\Einheit
\ifnum\@xarg =0 \@vvector
  \else \ifnum\@yarg =0 \@hvector \else \@svector\fi
\fi
\else\@badlinearg\fi}

\def\@hvector{\@hline\hbox to 0pt{\@linefnt
\ifnum \@xarg <0 \@getlarrow(1,0)\hss\else
    \hss\@getrarrow(1,0)\fi}}

\def\@vvector{\ifnum \@yarg <0 \@downvector \else \@upvector \fi}

\def\@svector{\@sline
\@tempcnta\@yarg \ifnum\@tempcnta <0 \@tempcnta=-\@tempcnta\fi
\ifnum\@tempcnta <5
  \hskip -\wd\@linechar
  \@upordown\@clnht \hbox{\@linefnt  \if@negarg
  \@getlarrow(\@xarg,\@yyarg) \else \@getrarrow(\@xarg,\@yyarg) \fi}%
\else\@badlinearg\fi}

\def\@upline{\hbox to \z@{\hskip -.5\Pfadd@cke \vrule width \Pfadd@cke
   height \@linelen depth \z@\hss}}

\def\@downline{\hbox to \z@{\hskip -.5\Pfadd@cke \vrule width \Pfadd@cke
   height \z@ depth \@linelen \hss}}

\def\@upvector{\@upline\setbox\@tempboxa\hbox{\@linefnt\char'66}\raise
     \@linelen \hbox to\z@{\lower \ht\@tempboxa\box\@tempboxa\hss}}

\def\@downvector{\@downline\lower \@linelen
      \hbox to \z@{\@linefnt\char'77\hss}}

\def\@getlarrow(#1,#2){\ifnum #2 =\z@ \@tempcnta='33\else
\@tempcnta=#1\relax\multiply\@tempcnta \sixt@@n \advance\@tempcnta
-9 \@tempcntb=#2\relax\multiply\@tempcntb \tw@
\ifnum \@tempcntb >0 \advance\@tempcnta \@tempcntb\relax
\else\advance\@tempcnta -\@tempcntb\advance\@tempcnta 64
\fi\fi\char\@tempcnta}

\def\@getrarrow(#1,#2){\@tempcntb=#2\relax
\ifnum\@tempcntb < 0 \@tempcntb=-\@tempcntb\relax\fi
\ifcase \@tempcntb\relax \@tempcnta='55 \or
\ifnum #1<3 \@tempcnta=#1\relax\multiply\@tempcnta
24 \advance\@tempcnta -6 \else \ifnum #1=3 \@tempcnta=49
\else\@tempcnta=58 \fi\fi\or
\ifnum #1<3 \@tempcnta=#1\relax\multiply\@tempcnta
24 \advance\@tempcnta -3 \else \@tempcnta=51\fi\or
\@tempcnta=#1\relax\multiply\@tempcnta
\sixt@@n \advance\@tempcnta -\tw@ \else
\@tempcnta=#1\relax\multiply\@tempcnta
\sixt@@n \advance\@tempcnta 7 \fi\ifnum #2<0 \advance\@tempcnta 64 \fi
\char\@tempcnta}
%%%%%%%%%%%%%%%%%%%%%%%%%%%%%%%%%%%%%%%%%%%%%%%%%%%%%%%%%%%%%

\def\Diagonale(#1,#2)#3{\unskip\leavevmode
  \xcoord#1\relax \ycoord#2\relax
      \raise\ycoord \Einheit\hbox to0pt{\hskip\xcoord \Einheit
         \Line@(1,1){#3}\hss}}
\def\AntiDiagonale(#1,#2)#3{\unskip\leavevmode
  \xcoord#1\relax \ycoord#2\relax %\advance\xcoord by -0.05\relax
      \raise\ycoord \Einheit\hbox to0pt{\hskip\xcoord \Einheit
         \Line@(1,-1){#3}\hss}}
\def\Pfad(#1,#2),#3\endPfad{\unskip\leavevmode
  \xcoord#1 \ycoord#2 \thicklines\ZeichnePfad#3\endPfad\thinlines}
\def\ZeichnePfad#1{\ifx#1\endPfad\let\next\relax
  \else\let\next\ZeichnePfad
    \ifnum#1=1
      \raise\ycoord \Einheit\hbox to0pt{\hskip\xcoord \Einheit
         \vrule height\Pfadd@cke width1 \Einheit depth\Pfadd@cke\hss}%
      \advance\xcoord by 1
    \else\ifnum#1=2
      \raise\ycoord \Einheit\hbox to0pt{\hskip\xcoord \Einheit
        \hbox{\hskip-\PfadD@cke\vrule height1 \Einheit width\PfadD@cke depth0pt}\hss}%
      \advance\ycoord by 1
    \else\ifnum#1=3
      \raise\ycoord \Einheit\hbox to0pt{\hskip\xcoord \Einheit
         \Line@(1,1){1}\hss}
      \advance\xcoord by 1
      \advance\ycoord by 1
    \else\ifnum#1=4
      \raise\ycoord \Einheit\hbox to0pt{\hskip\xcoord \Einheit
         \Line@(1,-1){1}\hss}
      \advance\xcoord by 1
      \advance\ycoord by -1
    \fi\fi\fi\fi
  \fi\next}
\def\hSSchritt{\leavevmode\raise-.4pt\hbox to0pt{\hss.\hss}\hskip.2\Einheit
  \raise-.4pt\hbox to0pt{\hss.\hss}\hskip.2\Einheit
  \raise-.4pt\hbox to0pt{\hss.\hss}\hskip.2\Einheit
  \raise-.4pt\hbox to0pt{\hss.\hss}\hskip.2\Einheit
  \raise-.4pt\hbox to0pt{\hss.\hss}\hskip.2\Einheit}
\def\vSSchritt{\vbox{\baselineskip.2\Einheit\lineskiplimit0pt
\hbox{.}\hbox{.}\hbox{.}\hbox{.}\hbox{.}}}
\def\DSSchritt{\leavevmode\raise-.4pt\hbox to0pt{%
  \hbox to0pt{\hss.\hss}\hskip.2\Einheit
  \raise.2\Einheit\hbox to0pt{\hss.\hss}\hskip.2\Einheit
  \raise.4\Einheit\hbox to0pt{\hss.\hss}\hskip.2\Einheit
  \raise.6\Einheit\hbox to0pt{\hss.\hss}\hskip.2\Einheit
  \raise.8\Einheit\hbox to0pt{\hss.\hss}\hss}}
\def\dSSchritt{\leavevmode\raise-.4pt\hbox to0pt{%
  \hbox to0pt{\hss.\hss}\hskip.2\Einheit
  \raise-.2\Einheit\hbox to0pt{\hss.\hss}\hskip.2\Einheit
  \raise-.4\Einheit\hbox to0pt{\hss.\hss}\hskip.2\Einheit
  \raise-.6\Einheit\hbox to0pt{\hss.\hss}\hskip.2\Einheit
  \raise-.8\Einheit\hbox to0pt{\hss.\hss}\hss}}
\def\SPfad(#1,#2),#3\endSPfad{\unskip\leavevmode
  \xcoord#1 \ycoord#2 \ZeichneSPfad#3\endSPfad}
\def\ZeichneSPfad#1{\ifx#1\endSPfad\let\next\relax
  \else\let\next\ZeichneSPfad
    \ifnum#1=1
      \raise\ycoord \Einheit\hbox to0pt{\hskip\xcoord \Einheit
         \hSSchritt\hss}%
      \advance\xcoord by 1
    \else\ifnum#1=2
      \raise\ycoord \Einheit\hbox to0pt{\hskip\xcoord \Einheit
        \hbox{\hskip-2pt \vSSchritt}\hss}%
      \advance\ycoord by 1
    \else\ifnum#1=3
      \raise\ycoord \Einheit\hbox to0pt{\hskip\xcoord \Einheit
         \DSSchritt\hss}
      \advance\xcoord by 1
      \advance\ycoord by 1
    \else\ifnum#1=4
      \raise\ycoord \Einheit\hbox to0pt{\hskip\xcoord \Einheit
         \dSSchritt\hss}
      \advance\xcoord by 1
      \advance\ycoord by -1
    \fi\fi\fi\fi
  \fi\next}
\def\Koordinatenachsen(#1,#2){\unskip
 \hbox to0pt{\hskip-.5pt\vrule height#2 \Einheit width.5pt depth1 \Einheit}%
 \hbox to0pt{\hskip-1 \Einheit \xcoord#1 \advance\xcoord by1
    \vrule height0.25pt width\xcoord \Einheit depth0.25pt\hss}}
\def\Koordinatenachsen(#1,#2)(#3,#4){\unskip
 \hbox to0pt{\hskip-.5pt \ycoord-#4 \advance\ycoord by1
    \vrule height#2 \Einheit width.5pt depth\ycoord \Einheit}%
 \hbox to0pt{\hskip-1 \Einheit \hskip#3\Einheit 
    \xcoord#1 \advance\xcoord by1 \advance\xcoord by-#3 
    \vrule height0.25pt width\xcoord \Einheit depth0.25pt\hss}}
\def\Gitter(#1,#2){\unskip \xcoord0 \ycoord0 \leavevmode
  \LOOP\ifnum\ycoord<#2
    \loop\ifnum\xcoord<#1
      \raise\ycoord \Einheit\hbox to0pt{\hskip\xcoord \Einheit\Punkt\hss}%
      \advance\xcoord by1
    \repeat
    \xcoord0
    \advance\ycoord by1
  \REPEAT}
\def\Gitter(#1,#2)(#3,#4){\unskip \xcoord#3 \ycoord#4 \leavevmode
  \LOOP\ifnum\ycoord<#2
    \loop\ifnum\xcoord<#1
      \raise\ycoord \Einheit\hbox to0pt{\hskip\xcoord \Einheit\Punkt\hss}%
      \advance\xcoord by1
    \repeat
    \xcoord#3
    \advance\ycoord by1
  \REPEAT}
\def\Label#1#2(#3,#4){\unskip \xdim#3 \Einheit \ydim#4 \Einheit
  \def\lo{\advance\xdim by-.5 \Einheit \advance\ydim by.5 \Einheit}%
  \def\llo{\advance\xdim by-.25cm \advance\ydim by.5 \Einheit}%
  \def\loo{\advance\xdim by-.5 \Einheit \advance\ydim by.25cm}%
  \def\o{\advance\ydim by.25cm}%
  \def\ro{\advance\xdim by.5 \Einheit \advance\ydim by.5 \Einheit}%
  \def\rro{\advance\xdim by.25cm \advance\ydim by.5 \Einheit}%
  \def\roo{\advance\xdim by.5 \Einheit \advance\ydim by.25cm}%
  \def\l{\advance\xdim by-.30cm}%
  \def\r{\advance\xdim by.30cm}%
  \def\lu{\advance\xdim by-.5 \Einheit \advance\ydim by-.6 \Einheit}%
  \def\llu{\advance\xdim by-.25cm \advance\ydim by-.6 \Einheit}%
  \def\luu{\advance\xdim by-.5 \Einheit \advance\ydim by-.30cm}%
  \def\u{\advance\ydim by-.30cm}%
  \def\ru{\advance\xdim by.5 \Einheit \advance\ydim by-.6 \Einheit}%
  \def\rru{\advance\xdim by.25cm \advance\ydim by-.6 \Einheit}%
  \def\ruu{\advance\xdim by.5 \Einheit \advance\ydim by-.30cm}%
  #1\raise\ydim\hbox to0pt{\hskip\xdim
     \vbox to0pt{\vss\hbox to0pt{\hss$#2$\hss}\vss}\hss}%
}
\catcode`\@=13

\catcode`\@=11
\font@\twelverm=cmr10 scaled\magstep1
\font@\twelveit=cmti10 scaled\magstep1
\font@\twelvebf=cmbx10 scaled\magstep1
\font@\twelvei=cmmi10 scaled\magstep1
\font@\twelvesy=cmsy10 scaled\magstep1
\font@\twelveex=cmex10 scaled\magstep1

\newtoks\twelvepoint@
\def\twelvepoint{\normalbaselineskip15\p@
 \abovedisplayskip15\p@ plus3.6\p@ minus10.8\p@
 \belowdisplayskip\abovedisplayskip
 \abovedisplayshortskip\z@ plus3.6\p@
 \belowdisplayshortskip8.4\p@ plus3.6\p@ minus4.8\p@
 \textonlyfont@\rm\twelverm \textonlyfont@\it\twelveit
 \textonlyfont@\sl\twelvesl \textonlyfont@\bf\twelvebf
 \textonlyfont@\smc\twelvesmc \textonlyfont@\tt\twelvett
%Erg\"anzung des fetten Small-Capitals-Fonts:
%
 \ifsyntax@ \def\big##1{{\hbox{$\left##1\right.$}}}%
  \let\Big\big \let\bigg\big \let\Bigg\big
 \else
  \textfont\z@=\twelverm  \scriptfont\z@=\tenrm  \scriptscriptfont\z@=\sevenrm
  \textfont\@ne=\twelvei  \scriptfont\@ne=\teni  \scriptscriptfont\@ne=\seveni
  \textfont\tw@=\twelvesy \scriptfont\tw@=\tensy \scriptscriptfont\tw@=\sevensy
  \textfont\thr@@=\twelveex \scriptfont\thr@@=\tenex
        \scriptscriptfont\thr@@=\tenex
  \textfont\itfam=\twelveit \scriptfont\itfam=\tenit
        \scriptscriptfont\itfam=\tenit
  \textfont\bffam=\twelvebf \scriptfont\bffam=\tenbf
        \scriptscriptfont\bffam=\sevenbf
  \setbox\strutbox\hbox{\vrule height10.2\p@ depth4.2\p@ width\z@}%
  \setbox\strutbox@\hbox{\lower.6\normallineskiplimit\vbox{%
        \kern-\normallineskiplimit\copy\strutbox}}%
 \setbox\z@\vbox{\hbox{$($}\kern\z@}\bigsize@=1.4\ht\z@
 \fi
 \normalbaselines\rm\ex@.2326ex\jot3.6\ex@\the\twelvepoint@}

\font@\fourteenrm=cmr10 scaled\magstep2
\font@\fourteenit=cmti10 scaled\magstep2
\font@\fourteensl=cmsl10 scaled\magstep2
\font@\fourteensmc=cmcsc10 scaled\magstep2
\font@\fourteentt=cmtt10 scaled\magstep2
\font@\fourteenbf=cmbx10 scaled\magstep2
\font@\fourteeni=cmmi10 scaled\magstep2
\font@\fourteensy=cmsy10 scaled\magstep2
\font@\fourteenex=cmex10 scaled\magstep2
\font@\fourteenmsa=msam10 scaled\magstep2
\font@\fourteeneufm=eufm10 scaled\magstep2
\font@\fourteenmsb=msbm10 scaled\magstep2
\newtoks\fourteenpoint@
\def\fourteenpoint{\normalbaselineskip15\p@
 \abovedisplayskip18\p@ plus4.3\p@ minus12.9\p@
 \belowdisplayskip\abovedisplayskip
 \abovedisplayshortskip\z@ plus4.3\p@
 \belowdisplayshortskip10.1\p@ plus4.3\p@ minus5.8\p@
 \textonlyfont@\rm\fourteenrm \textonlyfont@\it\fourteenit
 \textonlyfont@\sl\fourteensl \textonlyfont@\bf\fourteenbf
 \textonlyfont@\smc\fourteensmc \textonlyfont@\tt\fourteentt
%Erg\"anzung des fetten Small-Capitals-Fonts:
%
 \ifsyntax@ \def\big##1{{\hbox{$\left##1\right.$}}}%
  \let\Big\big \let\bigg\big \let\Bigg\big
 \else
  \textfont\z@=\fourteenrm  \scriptfont\z@=\twelverm  \scriptscriptfont\z@=\tenrm
  \textfont\@ne=\fourteeni  \scriptfont\@ne=\twelvei  \scriptscriptfont\@ne=\teni
  \textfont\tw@=\fourteensy \scriptfont\tw@=\twelvesy \scriptscriptfont\tw@=\tensy
  \textfont\thr@@=\fourteenex \scriptfont\thr@@=\twelveex
        \scriptscriptfont\thr@@=\twelveex
  \textfont\itfam=\fourteenit \scriptfont\itfam=\twelveit
        \scriptscriptfont\itfam=\twelveit
  \textfont\bffam=\fourteenbf \scriptfont\bffam=\twelvebf
        \scriptscriptfont\bffam=\tenbf
  \setbox\strutbox\hbox{\vrule height12.2\p@ depth5\p@ width\z@}%
  \setbox\strutbox@\hbox{\lower.72\normallineskiplimit\vbox{%
        \kern-\normallineskiplimit\copy\strutbox}}%
 \setbox\z@\vbox{\hbox{$($}\kern\z@}\bigsize@=1.7\ht\z@
 \fi
 \normalbaselines\rm\ex@.2326ex\jot4.3\ex@\the\fourteenpoint@}

\font@\seventeenrm=cmr10 scaled\magstep3
\font@\seventeenit=cmti10 scaled\magstep3
\font@\seventeensl=cmsl10 scaled\magstep3
\font@\seventeensmc=cmcsc10 scaled\magstep3
\font@\seventeentt=cmtt10 scaled\magstep3
\font@\seventeenbf=cmbx10 scaled\magstep3
\font@\seventeeni=cmmi10 scaled\magstep3
\font@\seventeensy=cmsy10 scaled\magstep3
\font@\seventeenex=cmex10 scaled\magstep3
\font@\seventeenmsa=msam10 scaled\magstep3
\font@\seventeeneufm=eufm10 scaled\magstep3
\font@\seventeenmsb=msbm10 scaled\magstep3
\newtoks\seventeenpoint@
\def\seventeenpoint{\normalbaselineskip18\p@
 \abovedisplayskip21.6\p@ plus5.2\p@ minus15.4\p@
 \belowdisplayskip\abovedisplayskip
 \abovedisplayshortskip\z@ plus5.2\p@
 \belowdisplayshortskip12.1\p@ plus5.2\p@ minus7\p@
 \textonlyfont@\rm\seventeenrm \textonlyfont@\it\seventeenit
 \textonlyfont@\sl\seventeensl \textonlyfont@\bf\seventeenbf
 \textonlyfont@\smc\seventeensmc \textonlyfont@\tt\seventeentt
%Erg\"anzung des fetten Small-Capitals-Fonts:
%
 \ifsyntax@ \def\big##1{{\hbox{$\left##1\right.$}}}%
  \let\Big\big \let\bigg\big \let\Bigg\big
 \else
  \textfont\z@=\seventeenrm  \scriptfont\z@=\fourteenrm  \scriptscriptfont\z@=\twelverm
  \textfont\@ne=\seventeeni  \scriptfont\@ne=\fourteeni  \scriptscriptfont\@ne=\twelvei
  \textfont\tw@=\seventeensy \scriptfont\tw@=\fourteensy \scriptscriptfont\tw@=\twelvesy
  \textfont\thr@@=\seventeenex \scriptfont\thr@@=\fourteenex
        \scriptscriptfont\thr@@=\fourteenex
  \textfont\itfam=\seventeenit \scriptfont\itfam=\fourteenit
        \scriptscriptfont\itfam=\fourteenit
  \textfont\bffam=\seventeenbf \scriptfont\bffam=\fourteenbf
        \scriptscriptfont\bffam=\twelvebf
  \setbox\strutbox\hbox{\vrule height14.6\p@ depth6\p@ width\z@}%
  \setbox\strutbox@\hbox{\lower.86\normallineskiplimit\vbox{%
        \kern-\normallineskiplimit\copy\strutbox}}%
 \setbox\z@\vbox{\hbox{$($}\kern\z@}\bigsize@=2\ht\z@
 \fi
 \normalbaselines\rm\ex@.2326ex\jot5.2\ex@\the\seventeenpoint@}

\catcode`\@=13

\topmatter
\author  George E. Andrews\footnote{\hbox{Partially supported by 
National Science Foundation Grant  DMS 9870060.}}\\
Christian Krattenthaler\footnote{\hbox{Partially supported by the Austrian
Science Foundation FWF, grant P13190-MAT.}}\\
Luigi Orsina\\ Paolo Papi\endauthor
\bigskip
\address
\vskip 0.pt Department of Mathematics
\vskip 0.pt The Pennsylvania State University 
\vskip 0.pt 215 McAllister Building
\vskip 0.pt University Park PA 16802, USA  
\vskip 0.pt e-mail:{\rm \ andrews\@math.psu.edu} 
\vskip 5pt
\vskip 0.pt Institut f\"ur Mathematik der Universit\"at Wien,
\vskip 0.pt Strudlhofgasse 4, A-1090 Wien, Austria.
\vskip 0.pt e-mail:{KRATT\@Ap.Univie.Ac.At}
\vskip 3pt
\vskip 5.pt Dipartimento di Matematica, Istituto G. Castelnuovo
\vskip 0.pt Universit\`a di Roma ``La Sapienza"  
\vskip 0.pt Piazzale Aldo Moro 5
\vskip 0.pt 00185 Rome --- ITALY 
\vskip 0.pt e-mail:{\rm \ orsina\@mat.uniroma1.it,\hskip3 pt papi\@mat.uniroma1.it}
 \endaddress
\abstract We study the combinatorics of $ad$-nilpotent ideals of  a
Borel subalgebra of 
$sl(n+1,\Bbb C)$. We provide an inductive method for calculating the
class of nilpotence of these ideals and   formulas for the 
number of 
 ideals  having a given class of nilpotence. We study the
relationships between these results and the 
combinatorics of Dyck paths, based upon  a remarkable bijection
between $ad$-nilpotent ideals and Dyck paths. Finally, we 
propose a $(q,t)$-analogue of the Catalan number
$C_n$. These
$(q,t)$-Catalan numbers count on the one hand $ad$-nilpotent ideals
with respect to dimension and class of nilpotence, and on the other hand  
admit interpretations in terms of natural statistics on Dyck paths.
\endabstract 

\leftheadtext {G. E. Andrews, C. Krattenthaler, L. Orsina and P. Papi}
\rightheadtext {$ad$-nilpotent
$\frak b$-ideals in $sl(n)$}

\title  $ad$-nilpotent
$\frak b$-ideals in $sl(n)$ having a fixed class of nilpotence: combinatorics and enumeration 
\endtitle
\endtopmatter
\nologo
\def\endemo{\qed\enddemo}

\def\g{\frak g}
\def\h{\frak h}
\def\n{\frak n}
\def\bb{\frak b}
\def\D{\Delta}
\def\l{\lambda}
\def\Dp{\Delta^+}

\def\Da{\widehat\Delta}
\def\Dap{\widehat\Delta^+}
\def\d{\delta}
\def\r{R(\h)}

\def\a{\alpha}
\def\b{\beta}

\def\I{\Cal I^n}
\def\l{\lambda}
\def\d{\delta}

\def\i{{\frak i}}

\def\N{N(w)}

\def\gg{\gamma}

\def\p{\Phi}

\TagsOnRight

\def\formel{1.1}
\def\tabalgo{3.1}
\def\fastalgo{3.2}
\def\lM{3.3}
\def\lm{3.4}
\def\summe{4.1}
\def\decomp{4.2}
\def\Summe{4.3}
\def\formela{4.4}
\def\formelb{4.5}
\def\formelc{4.6}
\def\qtCat{6.1}
\def\eqdmin{6.2}
\def\eqdmax{6.3}
\def\eqDmin{6.4}
\def\eqDmax{6.5}
\def\ineqm{6.6}
\def\eqintm{6.7}
\def\ineqM{6.8}
\def\eqstufen{6.9}
\def\eqintM{6.10}
\def\eqintmin{6.11}
\def\ineqA{6.12}
\def\ineqa{6.13}
\def\ineqb{6.14}

\def\l{\lambda}
\def\v#1{\left\vert#1\right\vert}
\def\fl#1{\left\lfloor#1\right\rfloor}
\def\cl#1{\left\lceil#1\right\rceil}
\def\flkl#1{\lfloor#1\rfloor}
\def\clkl#1{\lceil#1\rceil}
\def\dmin#1#2{{\theta^{\text{min}}_{#1}(#2)}}
\def\dmax#1#2{{\theta^{\text{max}}_{#1}(#2)}}
\def\Dmin#1#2{{\Theta^{\text{min}}_{#1}(#2)}}
\def\Dmax#1#2{{\Theta^{\text{max}}_{#1}(#2)}}
\def\({\left(}
\def\){\right)}

\document
\heading\S1 Introduction\endheading

Let $\g$ be a  complex simple Lie algebra of rank $n$. Let
$\h\subset\g$ be a fixed Cartan subalgebra,  
$\D$ the corresponding  root system of $\g$. Fix a positive system
$\Dp$ in $\D$, and let  
$\Pi=\{\a_1,\ldots,\a_n\}$ be the corresponding basis of simple roots.
For each $\a\in \Dp$ let ${\frak g}_\a$ be the root space of $\g$ relative to $\a$, 
$\n = \bigoplus \limits_{\a\in \Dp}{\frak g}_\a$,  and $\bb$ be the Borel subalgebra
$\bb=\h\oplus \n$.\par
Let $\I$ denote the set of $ad$-nilpotent ideals (i.e., consisting of
$ad$-nilpotent elements) of $\bb$. 
 %Then  $\i\in\I$  if and only if $\i\subseteq \n$.
These ideals, together with the subclass $\I_{ab}$ of Abelian ideals, have been studied in \cite{6}. In that paper,  Kostant 
stated a
useful equivalence criterion for certain  decomposably-generated simple $K$-submodules of
$\Lambda(\frak g)$ in terms of $\I$ (here,
$K$ is a compact semi-simple Lie group, and
$\frak g$ is the complexification of its Lie algebra). Moreover, 
he used the set of Abelian ideals to describe the eigenspace
relative to the maximal eigenvalue of the Laplace-Beltrami operator of $K$.\par  The subclass
$\I_{ab}$  has been studied much more recently in 
\cite{7} in connection with  discrete series representations.
The latter paper was partly motivated by a striking enumerative result
due to D. Peterson: the Abelian ideals are $2^n$ in 
number, independently of the type of
$\frak g$. (In contrast, the cardinality of
$\I$ depends on the type; see \cite{1, Th.~3.1}). Even more surprising
is the proof of Peterson's result, which involves the affine 
Weyl group $\widehat W$ of $\g$. In \cite{1}, the  encoding of the
ideals through certain elements of the affine Weyl group 
has been generalized from
$\I_{ab}$ to the entire set $\I$ of $ad$-nilpotent ideals. There it was shown
that any such ideal determines in a 
combinatorial way  the set of ``inversions" of a unique element  
in $\widehat W$. 

The combinatorial methods used in \cite{1} entailed
the problem of enumerating the ideals of $\I$ with respect to 
class of nilpotence. By definition, the class of nilpotence of an
ideal $\i$, which we denote here by $n(\i)$, is the smallest number
$m$ such that $m$-fold bracketing of $\i$ with itself gives the zero ideal. 
(Thus, the Abelian ideals are exactly those with class of nilpotence
at most 1.)
A solution to the problem, 
for $\frak g$ of type $A_n$, was obtained in \cite{10},
where it was shown that the number of ideals in $\I$ with class of
nilpotence $k$ is given by  
$$
\sum_{0=i_0<i_1<\dots<i_k<i_{k+1}=n+1}\prod_{j=0}^{k-1}
\binom{i_{j+2}-i_j-1}{i_{j+1}-i_j}.\tag\formel$$

The purpose of this paper is to deepen and enhance the understanding
of the enumerative properties of $ad$-nilpotent ideals of a fixed Borel
subalgebra of $sl(n+1,\Bbb C)$. First of all, after having recalled
the algebraic preliminaries in Section~2, we
describe in Section~3 a fast 
combinatorial algorithm for the computation of the class of nilpotence
of a given ideal (see Proposition~3.2). (We remark that it is based on
a ``slow" algorithm, see Proposition~3.1, which has interesting
relations to the elements of the affine symmetric group that are
obtained by the main result of \cite{1} mentioned earlier; see the
remarks at the end of Section~3.)
 
This algorithm implies naturally a partition into subintervals of the
interval $[(0,\ldots,0),(n,n-1,\ldots,1)]$ in the Young lattice, see
Proposition~4.1 in Section~4. From
this partition, the formula (\formel) follows immediately, thus
providing a proof different from the one in \cite{1} (see
Theorem~4.2).

However, formula (\formel) gives much more. Since
expressions like the one in (\formel) occur in the theory of Dyck
paths, it links the enumeration of $ad$-nilpotent ideals to the
enumeration of Dyck paths. To be precise, we prove that there are as
many $ad$-nilpotent ideals of a fixed subalgebra of $sl(n+1,\Bbb C)$
with class of nilpotence $k$ as there are Dyck paths of length $2n+2$
with height $k+1$ (see Theorem~4.4). 
As there are numerous formulas available for the
number of these Dyck paths, we obtain immediately alternative
expressions for the number of these $ad$-nilpotent ideals, see
Theorem~4.5 in Section~4. In particular, formula (\formelc) 
must be preferred over
formula (\formel), as it is much simpler and computationally superior.
Curious outcomes of these results are, for example, the observation 
that the number of $ad$-nilpotent ideals with class of nilpotence at
most 2 (instead of 1, as in Peterson's result) is a Fibonacci number,
as well as the observation 
that the number of $ad$-nilpotent ideals with class of
nilpotence at most 3 is essentially a power of 3,
see Corollary~4.7.

In Section~5 we make the connection between $ad$-nilpotent ideals and
Dyck paths completely explicit, by exhibiting a bijection between
$ad$-nilpotent ideals in $sl(n+1,\Bbb C)$
with class of nilpotence $k$ and Dyck paths of length $2n+2$
with height $k+1$.

The subject of Section~6 is an apparently new $(q,t)$-analogue of
Catalan numbers. (In particular, it is unrelated to the
$(q,t)$-Catalan numbers of Garsia and Haiman \cite{4}.) 
It counts
$ad$-nilpotent ideals in $sl(n+1,\Bbb C)$ simultaneously 
with respect to dimension and class of nilpotence. As it results
directly from the earlier mentioned interval decomposition, it is
composed out of a rather straightforward $(q,t)$-extension of formula
(\formel), see Theorem~6.1. 
In terms of
combinatorics, for $q=1$ this $(q,t)$-Catalan number reduces to the
generating function for Dyck paths counted with respect to height,
whereas for $t=1$ it reduces to the generating function for Dyck paths
counted with respect to area. 

Our combinatorial analysis allows us to provide precise results
concerning the minimal and maximal dimension of an ideal with fixed
class of nilpotence, and the minimal and maximal class of nilpotence of an
ideal with fixed dimension. In terms of our $(q,t)$-analogue of the
Catalan number this amounts to determining the minimal and maximal
degree in the variable $q$ once the degree in $t$ is fixed, and vice
versa. All this is also found in Section~6, see Theorems~6.2 and 6.3.

\smallskip
We now fix the notation that will be used
throughout the paper. As usual, we denote the set of
integers by $\Bbb Z$.
For binomial coefficients,
we will use the following convention: Given integers
$m$ and $n$, we let
$$
\binom{m}{n} = 
\cases
\frac{m!}{(m-n)!\,n!} \quad & \text{if $m \geq n > 0$,} \\
1 \quad & \text{if $n=0$,}\\
0 \quad & \text{in any other case.} \\
\endcases
$$
Similarly, we define the $t$-binomial coefficient by
$$\bmatrix{m}\\{n}\endbmatrix_t=\cases
\frac{[m]!}{[m-n]!\,[n]!} \quad & \text{if $m \geq n > 0$,} \\
1 \quad & \text{if $n=0$,}\\
0 \quad & \text{in any other case,} \endcases
$$
where the $t$-factorial $[m]!$ is defined by 
$[m]!=[m][m-1]\cdots [1]$, $[0]!=1$, with $[i]=(t^i-1)/(t-1)$.

Finally, for a partition $\l=(\l_1,\ldots,\l_n)$,
$\l_1\geq\dots\geq\l_n\geq 0$, we write $|\l|$ for the {\it size} 
$\sum_{i=1}^n\l_i$ of $\l$. 
We will identify a partition $\l=(\l_1,\ldots,\l_n)$ with its {\it Ferrers
diagram}, which is the array of cells with $n$ left-justified rows, the
$i$th row being of length $\l_i$. For example, Figure~1 shows the
(Ferrers diagrams of the) 
partitions $(3,2,1)$ and $(3,1)$. The cell in the $i$th row and $j$th
column will be always identified with the pair $(i,j)$.

\midinsert
\vskip10pt
\vbox{
$$
\Pfad(0,0),222111\endPfad
\Pfad(0,0),121212\endPfad
\Pfad(0,1),1212\endPfad
\Pfad(0,2),12\endPfad
\hbox{\hskip3.5cm}
\Pfad(0,1),22111\endPfad
\Pfad(0,1),12112\endPfad
\Pfad(0,2),112\endPfad
\Pfad(1,2),2\endPfad
\Label\lo{\bullet}(1,1)
\Label\lo{\bullet}(3,2)
\hskip1.5cm$$
\centerline{\eightpoint Figure 1}
}
\vskip10pt
\endinsert

We call a cell of a diagram a corner cell, if there are no cells to
the right and to the bottom. For example, the corner cells of the
diagram (corresponding to the partition) $(3,1)$ are the cells labelled
$(1,3)$ and $(2,1)$ (which are marked by bullets in Figure~1).

\heading\S2 Algebraic preliminaries \endheading

Let $\i\in \I$, i.e., $\i$ is 
an $ad$-nilpotent ideal of our fixed Borel subalgebra $\bb$.
Clearly, we can write $\i$ as
$\i=\bigoplus\limits_{\a\in \p}{\frak g}_\a$, for some collection $\p$ of
positive roots. The
collection $\p$ encodes an ideal $\i\in\I$ if and only if
for all $\a\in\p$ and $\b\in\Dp$ such that $\a+\b$ is a root, we have
$\a+\b\in\p.$ If one endows $\Dp$ with the (restriction of the) usual
partial order on the root lattice, that is, for $\a,\b\in\Dp$ we let 
$\a\leq\b$ if and only if $ \b-\a=\sum_{\gg\in\Dp}c_\gg\gg$, for some
nonnegative integers $c_\gg$, then this can be phrased differently as
follows: 
$\p$ encodes an ideal if and only
if it is a dual order ideal in $\Dp$.

In the rest of the paper we will exclusively 
deal with $\g$ of type $A_n$, i.e., the Lie
algebra $sl(n+1,\Bbb C)$ of $(n+1)\times (n+1)$ traceless matrices. The last
observation of the previous paragraph allows us to represent
$ad$-nilpotent ideals conveniently in a geometric fashion, which will
be crucial in all subsequent considerations (see also \cite{1,
Sec.~3}). Clearly, any positive
root in $A_n$ can be written as a sum of simple roots. Explicitly,  
let us write
$\tau_{ij}=\a_i+\dots+\a_{n-j+1},$  $1\leq i\le n$, $1\le j\leq
n-i+1$. If we place the roots $\tau_{ij}$, $j=1,2,\dots,n-i+1$, in the
$i$th row of a diagram, then this defines an arrangement of
the positive roots in a staircase fashion. For example, for
$A_3$ we obtain the arrangement
$$\matrix\format\l\quad &\l\quad &\l\\
\a_1+\a_2+\a_3 &\a_1+\a_2&\a_1\\
\a_2+\a_3 &\a_2 & \\
\a_3 \endmatrix$$
Obviously, the above defines an identification of positive
roots with cells of the staircase diagram $(n,n-1,\dots,1)$, 
in which the root $\tau_{ij}$
is identified with the cell $(i,j)$. For example, for $A_3$, the root
$\a_1+\a_2$ is identified with cell $(1,2)$ in the diagram $(3,2,1)$,
shown on the left in Figure~1.

Given an $ad$-nilpotent ideal $\i$, written as 
$\i=\bigoplus\limits_{\a\in \p}{\frak g}_\a$, for some collection $\p$ of
positive roots, we can use the above
identification to represent $\i$ as the set of cells that correspond
to the roots in $\p$. Since, as we noted above, $\i$ is a dual order
ideal, the set of cells obtained forms a (Ferrers diagram of a)
partition. For example, the ideal ${\frak g}_{\a_1+\a_2+\a_3}\oplus
{\frak g}_{\a_1+\a_2}\oplus {\frak g}_{\a_1}\oplus {\frak g}_{\a_2+\a_3}$ corresponds to the
partition $(3,1,0)$, shown on the right in Figure~1. 
Clearly, this correspondence is reversible as long as the partition is
contained in $(n,n-1,\dots,1)$. Thus, we
have defined a bijection between $ad$-nilpotent ideals in $sl(n+1,\Bbb
C)$ and subdiagrams of $(n,n-1,\dots,1)$.
In particular, as it is well-known that the number of the latter
subdiagrams is $C_{n+1}$, $C_n=
\frac{1}{n+1}\binom{2n}{n}$ being the $n$-th Catalan 
number, the number of $ad$-nilpotent
ideals is equal to the $(n+1)$st Catalan number 
(see \cite{12, Sec.~2} and also \cite{1, Sec.~3}).

\heading\S3 Calculating the class of nilpotence\endheading

The goal of this section is to describe a fast algorithm to determine
the class of nilpotence for any given $ad$-nilpotent ideal $\i$. As a first
step, we describe a tableau algorithm which computes the descending
central series of $\i$ (i.e., the $m$-fold bracketings of $\i$ with
itself, for any $m$). More precisely, let $t_{i,j}$ be the maximal
number $m$ such that the root space ${\frak g}_{\tau_{ij}}$ occurs in 
$$\i^{m}:=\underbrace {[\cdots[\i,\i],\dots]}_{m\text{ occurrences of
}\i}.$$
Then we claim that the numbers $t_{i,j}$ can be obtained as follows.
Let $\l$ be the subdiagram of $(n,n-1,\dots,1)$ that corresponds to
$\i$ according to the identification explained in Section~2. Define a
filling $(t_{i,j})_{1\le i\le n,\,1\le j\le n-i+1}$ 
of the cells of $(n,n-1,\dots,1)$ by recursively setting
$$t_{i,j}=\cases 0\quad &\text{if\ } (i,j)\notin \lambda,\\1\quad &\text{if\
} (i,j)\text{ is a corner cell of }\lambda,\\
\max\limits_{j<k\leq n-i+1}\ \{t_{i,k}+t_{n-k+2,j}\}\quad
&\text{otherwise.}\endcases
\tag\tabalgo$$ 
It is easy to see that the above rule uniquely defines a filling of
$(n,n-1,\dots,1)$, whose nonzero entries are precisely those
corresponding to the cells of $\l$. E.g., when $n=4$, the fillings 
corresponding to 
$(2,1,0,0)$, $(3,3,2,1)$, $(4,3,2,1)$ are respectively
$$
\matrix
1&1&0&0\\1&0&0\\0&0\\0\endmatrix\qquad
\matrix
3&2&1&0\\3&2&1\\2&1\\1\endmatrix\qquad
\matrix
4&3&2&1\\3&2&1\\2&1\\1\endmatrix$$

For the verification of our claim it suffices to observe that if
${\frak g}_{\a}\subseteq \i^{a}$ and ${\frak g}_{\b}\subseteq \i^{b}$ then, under
the assumption that $\a+\b$ is a root,
${\frak g}_{\a+\b}\subseteq \i^{a+b}$. For, with our labelling of positive
roots, a sum $\tau_{i,k}+\tau_{l,j}$ is a root if and only if
$l=n-k+2$ (or $i=n-j+2$, which, however, is the same case by
symmetry), in which case we have
$\tau_{i,k}+\tau_{n-k+2,j}=\tau_{i,j}$. 
Thus, if, for some $k$ with $j<k\le n-i+1$, 
we know that ${\frak g}_{\tau_{i,k}}$ occurs in $\i^{t_{i,k}}$ and that
${\frak g}_{\tau_{n-k+2,j}}$ occurs in $\i^{t_{n-k+2,j}}$, then it follows that 
${\frak g}_{\tau_{i,j}}$ occurs in $\i^{t_{i,k}+t_{n-k+2,j}}$. Clearly, the
maximum of all possible numbers $t_{i,k}+t_{n-k+2,j}$, is
equal to the maximal possible exponent $m$ such that
${\frak g}_{\tau_{ij}}\subseteq \i^{m}$. This is exactly the content of
(\tabalgo). 

Let us summarize our findings so far in the proposition below.

\proclaim{Proposition 3.1} Let $\i\in\I$. Then, for any $(i,j)$ with
$1\le i\le n$ and $1\le j\le n-i+1$, the maximal number $t_{i,j}$ such
that ${\frak g}_{\tau_{ij}}\subseteq \i^{t_{i,j}}$ can be determined by the
tableau algorithm given in {\rm(\tabalgo)}. In particular, the class of
nilpotence of $\i$ is equal to $t_{1,1}$, the entry in the top-left
cell. 
\endproclaim
%\demo{Proof} It follows easily from the definitions that
% $$n(\i_\p)=\max_{\a\in\p}\left\{k_\a\mid \a=
%\b_1+\cdots+\b_{k_\a},\ \b_j\in\p\right\}.$$
% We are going to prove, by induction on $n-i-j+2$ that 
%$t_{ij}$ is the maximum number of summands in $\p$ in a decomposition
%of $\tau_{ij}$ as a sum of roots in $\Dp$. From this 
%fact, we obtain the thesis by using the fact that $(t_{ij}^{A_\p})$
%is  both column-strict and row-strict.\par 
% Assume $n-i-j+2=1$; then
%$(i,j)$ belongs to the $n$-th antidiagonal of $T_n$, hence either it
%does not belong to 
%$A$ or it belongs to $\partial A$. In the latter case $\tau_{ij}$ is simple, hence indecomposable. In any case, the
%assignments for
%$t_{ij}$ affords the desired number of summands. Now suppose $n-i-j+2>1$; observe, for $(i,j)\in A$,  that if $\tau_{ij}$ does
%not admit a decomposition as a sum of two roots at least one of which belongs to  $\p$, then $(i,j)\in \partial A$, and the
%definition of $t_{ij}$ gives $t_{ij}=1$,  as it should be. Otherwise, consider the decompositions of $\tau_{ij}$ as a sum of
%two roots: $\tau_{ij}=\tau_{ik}+\tau_{n-k+2\,j},\ j<k\leq n-i+1$. Using the definition of $t_{ij}$ in this case
%we are done by the inductive assumption.\endemo

In view of the second statement of Proposition~3.1, this tableau
algorithm provides an algorithm for the determination of the class of
nilpotence, which, however, is rather slow, as it involves the
determination of {\it all\/} the entries in the filling $(t_{i,j})$.
We will now show that, if one is only interested in the determination
of $t_{1,1}$ (which, by the second statement of Proposition~3.1 gives
exactly the class of nilpotence), then a considerable speedup can be
achieved. For a convenient statement of the result, we write, in abuse
of notation, $n(\l_1,\l_2,\ldots,\l_n)$ for $n(\i)$, given that the
partition corresponding to $\i$ according to the construction in
Section~2 is $(\l_1,\l_2,\ldots,\l_n)$.

\proclaim{Proposition 3.2}
Let $\i\in\I$ and let 
$\l=(\l_1,\ldots,\l_n)$ be the corresponding partition.
If $\l\ne(0,0,\ldots,0)$ then
$$n(\l_1,\l_2,\ldots,\l_n)=n(\l_{n+2-\l_1},\ldots,\l_n)+1.\tag\fastalgo$$
\endproclaim

It should be noted that on the left-hand side of (\fastalgo) appears
the class of nilpotence of an ideal in $\I$, whereas on the right-hand
side there appears the class of nilpotence of an ideal in $\Cal
I^{\l_1-1}$ (with corresponding partition
$(\l_{n+2-\l_1},\ldots,\l_n)$). 
The computation, however, can be carried out completely formally,
without reference to ideals, which we now demonstrate by an example.

\example{Example}
Let $\i\in\Cal I^{13}$ be the ideal which corresponds to the partition 
$(10,10,9,6,\mathbreak 5,4,4,3,1,1,1,1,0)$. (This is the partition in
Figure~2. At this point, all dotted lines should be ignored.) 
Then, by applying Proposition~3.2
iteratively, we obtain for the class of nilpotence of $\i$
$$\align
n(\i)&=n(10,10,9,6,5,4,4,3,1,1,1,1,0)\\ 
&= n(5,4,4,3,1,1,1,1,0)+1\\
&= n(1,1,1,0)+2\\
&=3.
\endalign$$
\endexample

As is obvious from the example, iterated application of
Proposition~3.2 provides a very efficient algorithm of determining the
class of nilpotence of a given ideal $\i$. 

Before we move on to the proof,
we wish to point out that this algorithm has a very nice geometric
rendering. Let, as before, $\l=(\l_1,\l_2,\dots,\l_n)$ 
be the partition corresponding to
$\i$. Consider the Ferrers diagram of $\l$. As it is contained in the
staircase diagram $(n,n-1,\dots,1)$, it must not cross the
antidiagonal line $x+y=n+1$.
We draw a zig-zag line
as follows (see Figure~2, where $n=13$ and
$\l=(10,10,9,6,5,4,4,3,1,1,1,1,0)$): 
we start on the vertical edge on the right of cell
$(1,\l_1)$, and move downward until we touch the antidiagonal
$x+y=n+1$. At the touching point we turn direction from 
vertical-down to horizontal-left, and move on until we touch a vertical
part of the Ferrers diagram. At the touching point we turn direction
from  horizontal-left to vertical-down. Now the procedure is
iterated, until we reach the line $x=0$. The class of nilpotence of
the ideal $\i$ is equal to the number of touching points on
$x+y=n+1$. In Figure~2, the resulting zig-zag line is the dotted line
outside the Ferrers diagram of $(10,10,9,6,5,4,4,3,1,1,1,1,0)$.
There are three touching points on $x+y=n+1=14$. (At this
point, the dotted lines inside the diagram should still be ignored.) 

\midinsert
\vskip10pt
\vbox{
$$
\Pfad(0,1),2222222222221111111111\endPfad
\Pfad(0,1),1222211212212121112122\endPfad
\SPfad(0,0),1222\endSPfad
\SPfad(1,4),11112222\endSPfad
\SPfad(5,9),11111222\endSPfad
\SPfad(0,0),2221\endSPfad
\SPfad(1,4),22221111\endSPfad
\SPfad(5,9),22211111\endSPfad
\thinlines
\Diagonale(0,-1){14}
\Label\r{x+y=n+1}(15,12)
\DuennPunkt(1,0)
\DuennPunkt(5,4)
\DuennPunkt(10,9)
\hskip7cm
$$
\centerline{\eightpoint Figure 2}
}
\vskip10pt
\endinsert

Since we need it in the proof of Proposition~3.2, and also later, 
let us express this
geometric rendering in formal terms. Obviously, the zig-zag line
describes the shape of a partition 
$$
\left( i_k^{n-i_k+1},
i_{k-1}^{i_k-i_{k-1}},\ldots,i_1^{i_2-i_1},0^{i_1-1}\right),
\tag\lM$$
where $i_k=\l_1$, $i_{k-1}=\l_{n-i_k+2}$, $i_{k-2}=\l_{n-i_{k-1}+2}$,
\dots, $i_1=\l_{n-i_2+2}$. Clearly we have $0<i_1<\dots<i_k<n+1$.
(In Figure~2, we have $k=3$ and $i_3=10$, $i_2=5$, $i_1=1$.)
Any partition $\l$ which gives rise to this zig-zag line must
necessarily contain the cells $(1,i_k)$, $(n-i_k+2,i_{k-1})$,
$(n-i_{k-1}+2,i_{k-2})$, \dots, $(n-i_2+2,i_1)$. (In Figure~2, these
are the cells $(1,10)$, $(5,5)$, $(10,1)$.) The ``minimal" partition
(in the sense of inclusion of diagrams) which contains these cells is 
$$
\left( i_k, i_{k-1}^{n-i_k+1},
i_{k-2}^{i_k-i_{k-1}},\ldots,i_1^{i_3-i_2},0^{i_2-2}\right). 
\tag\lm$$
(In Figure~2, this ``minimal" partition is indicated by the dotted
lines inside the Ferrers diagram of $(10,10,9,6,5,4,4,3,1,1,1,1,0)$.)
For later use, let us denote the partition in (\lM) by $\l_{i_1,\ldots,
i_k}^M$, and the partition in (\lm) by $\l_{i_1,\ldots, i_k}^m$.

\demo{Proof of Proposition 3.2}
We first show that the class of nilpotence is at least as large as the
number, $k$ say, which is typed out by the algorithm, in its geometric
rendering. 
Let $i_1,\dots,i_k$ be as above, $0<i_1<\dots<i_k<n+1$. As we
already noted, the partition $\l$ contains the cells 
$(1,i_k)$, $(n-i_k+2,i_{k-1})$,
$(n-i_{k-1}+2,i_{k-2})$, \dots, $(n-i_2+2,i_1)$. In view of the
correspondence of Section~2, these cells correspond to the root spaces 
${\frak g}_{\tau_{1,i_k}}$, ${\frak g}_{\tau_{n-i_k+2,i_{k-1}}}$, \dots,
${\frak g}_{\tau_{n-i_2+2,i_1}}$ contained in the ideal $\i$. The bracket of
$$[\cdots[{\frak g}_{\tau_{1,i_k}},{\frak g}_{\tau_{n-i_k+2,i_{k-1}}}],\dots,
{\frak g}_{\tau_{n-i_2+2,i_1}}]$$
is simply ${\frak g}_{\tau_{1,i_1}}$. (In particular, it is nontrivial.) 
Hence, the class of nilpotence of $\i$
is at least $k$, as was claimed.

In order to see that the class of nilpotence does not exceed $k$, we
consider the ideal, $\i^M$ say, 
which corresponds to the partition $\l_{i_1,\ldots,
i_k}^M$ (see (\lM)). Clearly, this ideal contains $\i$. Hence, its
class of nilpotence is an upper bound for the class of nilpotence of
$\i$. However, as is seen by inspection, the tableau algorithm
(\tabalgo) yields the following for $\l_{i_1,\ldots,
i_k}^M$: the entry $t_{i,j}$, where $n-i_{s}+2\le i<n-i_{s-1}+2$ and
$i_{r-1}<j\le i_{r}$, is given by
$s-r$. (Here, by convention, we have put $i_0:=0$ and $i_{k+1}:=n+1$.) 
In particular, the top-left entry, $t_{1,1}$, which by
Proposition~3.2 yields the class of nilpotence of $\i^M$, 
is equal to $(k+1)-1=k$. Hence, the class of nilpotence of $\i$ cannot
exceed $k$, and thus must be equal to $k$.
\endemo

\medskip
At the end of this section, we want to relate Proposition~3.1 to the
main result from \cite{1}, the latter setting up a connection between
$ad$-nilpotent ideals and elements of the affine Weyl group for any
type of $\frak g$.

Let $\Da$ and
$\widehat W$ be the affine real root system and the affine  Weyl group associated to $\D$ \cite{5}.
Having  fixed a positive system  $\Dp$ in $\D$, we have a corresponding positive system
$\Dap=(\Dp+\Bbb N\d)\cup(-\Dp+\Bbb Z^+\d)$  in $\Da$ ($\d$ is the "imaginary root").\par 
For $\p\subseteq\Dp$, set 
$\p^k=\left(\p^{k-1}+\p\right)\cap\D.$
If moreover $\p$ is a dual order ideal in $\Dp$ (cf\. the first
paragraph of Section~2),  define
$${\frak g}_{\p}=\bigcup\limits_{k\in\Bbb Z^{+}} (-\p^{k}+k\d).$$

Set 
$\N=\{\a\in\Dap\mid w^{-1}(\a)\in  \Da^-\},$
where $\Da^-=-\Dap$. It is well known that $N(w)$ determines $w$ uniquely. The main result of \cite{1} is the following
theorem, which holds for any simple Lie algebra. 
\proclaim{Theorem 3.3} Consider the ideal $\i_\p\in\I$ defined by
$\i_\p=\bigoplus\limits_{\a\in \p}{\frak g}_\a$, where $\p$ is a dual order
ideal in $\Dp$.
Then there exists a unique $w_\p\in\widehat W$
such that ${\frak g}_\p=N(w_\p)$. Moreover ${\frak g}_\p$ is the minimal set of the form $N(v),\,v\in\widehat W$
(w.r.t. inclusion) containing $-\p+\d$.\endproclaim

\proclaim{Proposition 3.4} Let $\i_\p\in\I$ be as in Theorem~{\rm3.3} and
let $(t_{i,j})$ be the corresponding filling of $(n,n-1,\dots,1)$
{\rm(}cf\. Section~{\rm3)}. Then
$$N(w_\p)=\bigcup_{1\leq i\leq j\leq n} \left\{-\tau_{ij}+h\d\mid
1\leq h\leq t_{i,j}\right\}.$$ 
\endproclaim
\demo{Proof} This is immediate from Proposition~3.1.\endemo 

This result admits the following interpretation. Recall that, in type
$\widetilde A_n$, $\widehat W$ can be realized as the 
group of {\it affine permutations} \cite{8}:
$$\align \widehat W\cong\bigg\{w:\Bbb Z\leftrightarrow\Bbb Z\mid
&w(t+n+1)=w(t)+n+1\ \forall\,t\in \Bbb Z,\quad\\ 
&\sum_{t=1}^{n+1}w(t)=\frac{(n+2)(n+1)}{2}\bigg\}.\endalign$$
Proposition 3.4, together with \cite{11, Th.~1} or \cite{13, Th.~3.2},
shows that the filling $(t_{i,j})$ determines the {\it 
inversion table}
\cite{2, Sec.~8} of 
$w_\p$, thought of as an affine permutation. In different terms,
$$\left\lfloor\frac{w_\p^{-1}(j)-w_\p^{-1}(i)}{n+1}\right\rfloor=
t_{i, n-j+2}, \qquad 1\leq i<j\leq n+1.$$
One can be even more explicit by using \cite{13, Th.~5.2}. Namely, 
for $1\leq i\leq n+1$, we have
$$w_\p^{-1}(i)=i+\sum_{j=1}^{i-1}t_{j,n-i+2}-\sum_{j=i+1}^{n+1}t_{i,n-j+2}.
$$
\proclaim{Problem} Find a combinatorial characterization of the affine
permutations $w_\p$ that correspond, as described above, 
to $ad$-nilpotent ideals.\endproclaim

\heading\S4 Enumeration of $ad$-nilpotent ideals\endheading

In this section we provide several formulae for the number of $ad$-nilpotent
ideals having a fixed class of nilpotence. The point of departure is 
a remarkable partition of the interval
$[(1,\ldots,0),(n,n-1,\ldots,1)]$ in the Young lattice that is implied
by the fast algorithm for the determination of
class of nilpotence given in Proposition~3.2.

\proclaim{Proposition 4.1}{\bf (a)} The interval
$I=[(0,\ldots,0),(n,n-1,\ldots,1)]$ in the Young lattice can be
decomposed into disjoint subintervals as  
$$I=\bigcup_{k=0}^n\,\,\bigcup_{0=i_0<i_1<\dots<i_k<i_{k+1}=n+1}
[\l^m_{i_1,\dots,i_k},\l^M_{i_1,\dots,i_k}],$$
where $\l^m_{i_1,\dots,i_k}$ and $\l^M_{i_1,\dots,i_k}$ are defined by
{\rm(\lm)} and {\rm(\lM)}, respectively. 

\noindent{\bf (b)} Let $\i$ be an ideal with corresponding partition
$\l$. If $\l\in[\l^m_{i_1,\dots,i_k},\l^M_{i_1,\dots,i_k}]$,
then the class of nilpotence of $\i$ is equal to $k$.\endproclaim 
\demo{Proof} 
This follows immediately from the geometric rendering of
Proposition~3.2 (see the remarks after the statement of
Proposition~3.2) and the arguments given in the proof of
Proposition~3.2.\endemo

\proclaim{Theorem 4.2} The number of ideals in $\I$ with class of
nilpotence $k$ is equal to
$$
\sum_{0=i_0<i_1<\dots<i_k<i_{k+1}=n+1}\prod_{j=0}^{k-1}\binom{i_{j+2}-i_j-1}{i_{j+1}-i_j}.
\tag\summe$$
\endproclaim
\demo{Proof} 
By Proposition 4.1 we have to count 
the number of partitions $\l$ with\linebreak
$\l_{i_1,\ldots, i_k}^m\subseteq
\l\subseteq \l_{i_1,\ldots, i_k}^M$, when $i_1,\ldots, i_k$ vary. 
For fixed $i_1,\ldots, i_k$ the corresponding
number is easily determined: the interval $[\l_{i_1,\ldots, i_k}^m,
\l_{i_1,\ldots, i_k}^M]$ decomposes into the product of $k$ Young
lattices as follows,
$$\multline
\left[\emptyset,\((i_k-i_{k-1})^{n-i_k}\)\right]\times
\left[\emptyset,\((i_{k-1}-i_{k-2})^{i_k-i_{k-1}-1}\)\right]\\
\times
\left[\emptyset,\((i_{k-2}-i_{k-3})^{i_{k-1}-i_{k-2}-1}\)\right]\times\dots
\times\left[\emptyset,\(i_1^{i_2-i_1-1}\)\right].
\endmultline\tag\decomp$$
(Here, $\emptyset$ stands for the empty partition.)
This decomposition is most obvious from Figure~2. There, the dotted lines
mark the partitions $\l_{i_1,\ldots, i_k}^m$ and $\l_{i_1,\ldots,
i_k}^M$ (with $k=3$, $i_3=10$, $i_2=5$, $i_1=1$). As is obvious from
the picture, the dotted lines determine $k$ (in Figure~2, we have $k=3$)
``independent" rectangles. So, if $\l\in[\l_{i_1,\ldots, i_k}^m,
\l_{i_1,\ldots, i_k}^M]$, there is only freedom within the rectangles,
which is expressed by the decomposition (\decomp). As the number of
partitions which are contained in a rectangle $(a^b)$ is equal to the
binomial coefficient $\binom {a+b}a$, the result follows.
\endemo

\proclaim{Corollary 4.3}  
The number of ideals in $\I$ with class of nilpotence at most $h$ is
equal to
$$
\sum_{0=i_0\le i_1\le\dots\le i_h\le
i_{h+1}=n+1}\prod_{j=0}^{h-1}\binom{i_{j+2}-i_j-1}{i_{j+1}-i_j}.
\tag\Summe$$
\endproclaim

\demo{Proof}
According to Theorem~4.2, we have to sum the expression (\summe) over
$k$ from $0$ to $h$. Because of our convention for binomial
coefficients (cf\. the introduction), this does indeed yield
(\Summe). For, if in (\Summe) we encounter $i_j$ and $i_{j+1}$ with
$i_j=i_{j+1}$ and $j\ge1$, then the binomial coefficient $\binom
{i_{j+1}-i_{j-1}-1}{i_j-i_{j-1}}$, which occurs in the summand,
vanishes. Hence, the only nonzero contributions in (\Summe) are by
indices $0=i_0= i_1=\dots=i_{h-k}<i_{h-k+1}<\dots< i_h<
i_{h+1}=n+1$, for some $k$. Because of our convention that $\binom
{-1}0=1$, the corresponding summand reduces to a term which appears in
the sum (\summe), upon replacement of $i_j$ by $i_{j-h+k}$,
$j=h-k,h-k+1,\dots,h+1$.
\endemo

This corollary makes the link of the enumeration of $ad$-nilpotent
ideals in $sl(n+1,\Bbb C)$ to the enumeration of Dyck paths. 
Recall that a Dyck path
is a lattice path  from $(0,0)$ to $(2n,0)$ with diagonal step vectors
$(1,1)$ and $(1,-1)$ which does not pass below the $x$-axis. 
We define the height
of a Dyck path to be the maximum ordinate of its peaks.

\proclaim{Theorem 4.4}
The number of ideals in $\I$ with class of nilpotence $k$ 
is exactly the same as the number of
Dyck paths from $(0,0)$ to $(2n+2,0)$ with height $k+1$.
\endproclaim
\demo{Proof}
The expression (\Summe) occurs in \cite{3, Proposition~3.B}. (There,
replace $n$ by $n+1$, $n_j$ by $i_{h-j+1}-i_{h-j}$,
$a_jb_{j+1}$ by $x$, $j=0,1,\dots,h$, and extract the coefficient of
$x^{n+1}$.) If this is combined with Corollary~2 in \cite{3}, then it
follows that the expression (\Summe) is equal to the number of Dyck
paths from $(0,0)$ to $(2n+2,0)$ with height at most $h+1$. Clearly,
since by Corollary~4.3 we know that it also equals the number of
ideals in $\I$ with class of nilpotence at most $h$, this implies the result.
\endemo

An immediate question is, of course, whether it is possible to provide
an explicit bijection between the ideals and Dyck paths in
Theorem~4.4. We are going to construct such a bijection in Section~5.
(It should be noted that the obvious correspondence between ideals and
partitions that we described in Section~2 cannot serve this
purpose. Although the border of 
a partition contained in $(n,n-1,\dots,1)$ can be
viewed as a Dyck path if the Ferrers diagram is rotated by $45^\circ$ in
the negative direction, this correspondence 
does not convert class of nilpotence of the ideal into height of
Dyck paths. For example, under this correspondence, the zero ideal,
the unique ideal with class of nilpotence $0$,
translates into the unique Dyck path with height $n+1$, i.e., the Dyck
path with $n+1$ up-steps followed by $n+1$ down-steps.)

\medskip
The enumeration of Dyck paths (and of lattice paths in general) is a
well-explored territory, where many explicit results exist. In view of
Theorem~4.4, these may now be used to obtain results for ideals with a
given class of nilpotence.

\proclaim{Theorem 4.5}
The number of ideals in $\I$ with class of nilpotence at most $h$ is
{\rm(}aside from {\rm(\Summe))}
equal to any of the following expressions:
$$\align
\det&\(\binom{i-\max\{0,j-h\}+1}{j-i+1}\)_{1\le i,j\le n}\tag\formela\\
&=\det\(\binom{i-j+h+1}{j-i+1}\)_{1\le i,j\le n}\tag\formelb\\
&= \sum_{k\in\Bbb Z} \frac {2k(h+3)+1} {2n+3}\binom{2n+3}
{n+1-k(h+3)}.\tag\formelc
\endalign$$ 
\endproclaim
\demo{Proof} 
We observe
that, instead of counting Dyck paths, we may equivalently
count lattice paths from $(0,0)$ to $(n+1,n+1)$ with step vectors
$(1,0),(0,1)$ which do not touch the lines 
$y=x-1$ and $y=x+h+2$. Then the determinantal expressions follow from 
\cite{9, Ch.~2, Th.~1}, while (\formelc) results from \cite{9, Ch.~1,
Th.~2} upon little simplification. \endemo 

Generating function results for Dyck paths translate into the
following result for ideals with a given class of nilpotence.

\proclaim{Theorem 4.6} Let $U_n(x)$ denote the $n$th Chebyshev
polynomial of the second 
kind, $U_n(\cos t)=\sin((n+1)t)/\sin t$, or, explicitly,
$$U_n(x)=\sum_{j\geq 0}(-1)^j\binom{n-j}{j}(2x)^{n-2j}.$$
Let $\a_n(h)$ denote the number of ideals in $\I$ with class of
nilpotence at most $h$. Then
$$\align
1+\sum_{n=0}^\infty\a_n(h) x^{n+1}&=
\frac{U_{h+1}\(1/2\sqrt x\)}{\sqrt x\,U_{h+2}\(1/2\sqrt x\)}\\
&=\cfrac 1\\1-
\cfrac x\\1-
\cfrac x\\
\cfrac \ddots\\1-x\endcfrac\endcfrac\endcfrac\endcfrac\ .
\endalign$$
(In the continued fraction there are $h+1$ occurrences of $x$.)
\endproclaim

\demo{Proof}
The expression in terms of a quotient of Chebyshev polynomials
follows, for example, from \cite{3, Prop.~12}, while the continued
fraction follows from Flajolet's continued fraction \cite{3, Th.~1}.
\endemo

By specializing this generating function result to $h=1$, we recover 
Peterson's result (in type $A_n$) that the number of Abelian
ideals (i.e., the ideals with class of nilpotence at most 1) in $\I$ is
$2^n$. If we specialize Theorem~4.6 to $h=2$ and $h=3$, we may obtain
further enumeration results, which are equally remarkable.

\proclaim{Corollary 4.7}
The number of $ad$-nilpotent ideals in $\I$ with class of nilpotence
at most $2$ is the Fibonacci number $F_{2n}$.
The number of $ad$-nilpotent ideals in $\I$ with class of nilpotence
at most $3$ is $(3^n+1)/2$. 
\endproclaim

\heading\S5 A bijection between $ad$-nilpotent ideals and Dyck
paths\endheading 

 Now we describe a bijection between ideals in $\I$ with class of
nilpotence $k$ and Dyck paths from $(0,0)$ to $(2n+2,0)$ with height
$k+1$. 

Let $\i\in\I$, and let $\l$ be the corresponding partition
contained in $(n,n-1,\dots,1)$, according to the correspondence
described in Section~2. The first step consists of determining
the interval, according to the decomposition of Proposition~4.1, the
partition $\l$ is in. I.e., we determine the integers
$i_k,i_{k-1},\dots,i_1$  such that
$\l\in[\l^m_{i_1,\dots,i_k},\l^M_{i_1,\dots,i_k}]$. 
To use the example of Section~3, 
$\l=(10,10,9,6,5,4,4,3,1,1,1,\mathbreak
1,0)$ (with $n=13$), we have $k=3$, $i_3=10$,
$i_2=5$, $i_1=1$. Figure~2 shows this partition. The dotted line
outside indicates the partition $\l^M_{i_1,\dots,i_k}=\l^M_{1,5,10}$,
the dotted line inside indicates $\l^m_{i_1,\dots,i_k}=\l^m_{1,5,10}$.

Now one generates a Dyck path step by step.
One starts with $n+1-i_k$ up-down pieces (in our example: $k=3$ and 
$n+1-i_k = n+1-i_3 = 4$; see Figure 3.a).

\midinsert
\vskip10pt
\vbox{
$$
\Einheit.36cm
\Gitter(9,3)(0,0)
\Koordinatenachsen(9,3)(0,0)
\Pfad(0,0),34343434\endPfad
\hskip2.88cm
$$
\centerline{\eightpoint a.}
$$
\Einheit.36cm
\Gitter(19,4)(0,0)
\Koordinatenachsen(19,4)(0,0)
\Pfad(0,0),343344334343443344\endPfad
\hskip6.48cm
$$
\centerline{\eightpoint b.}
$$
\Einheit.36cm
\Gitter(29,6)(0,0)
\Koordinatenachsen(29,6)(0,0)
\Pfad(0,0),3433344433433443334434443344\endPfad
\hskip10cm
$$
\centerline{\eightpoint c.}
\centerline{\eightpoint Figure 3}
}
\vskip10pt
\endinsert

In order to explain the next steps, we need to observe (as we did
already earlier) that the
interval
$[\l^m_{i_1,\dots,i_k},\l^M_{i_1,\dots,i_k}]$
pictorially decomposes into $k$ independent rectangles. In the example
of Figure~2, these are the rectangles formed by the dotted lines, the
top-most rectangle being a $3\times 5$ rectangle, the next a $4\times
4$ rectangle, and the bottom-most a $3\times 1$ rectangle (see Figure~2).

Now, in the top-most rectangle, we follow the shape of $\l$
inside the rectangle, from top-right to bottom-left. In our example of
Figure~2, this shape is $dldllldl$, the letter $d$ indicating a down-step
in the shape, the letter $l$ indicating a left-step. Thus, to the
portion of the shape contained in the rectangle corresponds a word
$l^{a_0}d\,l^{a_1}d\dots d\,l^{a_{i_k-i_{k-1}}}$. We insert $a_0$ up-down
pieces into the first peak of the already existing Dyck path (which,
by now, is just a zig-zag line; see Figure~3.a), 
$a_1$ up-down pieces into the second
peak, etc. In our example this generates the Dyck path in Figure~3.b.

This procedure is now repeated, by considering the remaining rectangles
one-by-one, from top to bottom. From now on, up-down pieces are only
inserted into {\it highest\/} peaks. To continue our example, the next
shape portion to be considered (the one contained in the $4\times 4$
rectangle) is $lddldlld$. Hence, $1$ up-down piece is inserted into the
first peak (of height $2$, since only highest peaks are
considered for insertions) 
in Figure~3.b, $0$ up-down pieces into the
second peak, etc.

The final result of this procedure, applied to the partition in
Figure~2, is shown in Figure~3.c (i.e., after also having considered
the bottom-most rectangle).

It is obvious that the result of this mapping is a Dyck path with
height $k+1$.
Conversely, given a Dyck path with height $k+1$, it is obvious how to
reverse the mapping and obtain the corresponding partition $\l$, and,
thus, the corresponding ideal $\i$ with
$n(\i)=k$. Therefore we have found the desired bijection.

\heading\S6 A  $(q,t)$-analogue of the Catalan number\endheading

As we said in the introduction, the total number of $ad$-nilpotent
ideals of a Borel subalgebra of $sl(n+1,\Bbb C)$ is the Catalan number
$C_{n+1}$. Let $\a_{n}(h,k)$ be the number of such ideals with
dimension $h$ and class
of nilpotence $k$. Then the generating function
$$C_n(q,t)=\sum_{h,k\geq 0}\a_n(h,k)t^hq^k$$
is a $(q,t)$-analogue of the Catalan number $C_{n+1}$. (It is unrelated
to the Garsia-Haiman $(q,t)$-Catalan number \cite{4}. This can be
seen, for example, by recalling that the
Garsia-Haiman $(q,t)$-Catalan number is symmetric in $q$ and $t$,
whereas our $(q,t)$-Catalan number is highly nonsymmetric.)

Define the {\it area} $A(P)$ of a Dyck path $P$ as the area of the region
between $P$ and the $x$-axis. Our $(q,t)$-Catalan number has the
following properties.

\proclaim{Theorem 6.1}
We have
$$C_n(q,t)=\sum_{k=0}^n
\left(\sum_{0=i_0<i_1<\dots<i_k<i_{k+1}=n+1} 
\prod_{j=0}^{k-1}t^{i_{j+1}(i_{j+3}-i_{j+2})}
\bmatrix {i_{j+2}-i_j-1}\\{i_{j+1}-i_j}\endbmatrix_t
\right)q^k,\tag\qtCat$$
with $i_{k+2}=n+2$. $C_n(q,1)$ is the generating function for
Dyck paths from $(0,0)$ to $(2n+2,0)$ counted with respect to height.
$C_n(1,t)$ is the
generating function for the same set of Dyck paths with respect to the
weight function ${(n+1)^2}/{2}-A(\cdot)$.
\endproclaim
\demo{Proof}
The expression (\qtCat) is obtained by following along the arguments
of the proof of Theorem~4.2. That is, for fixed class of nilpotence,
we use the decomposition (\decomp) (see also Figure~2). 
This reduces the problem to the
problem of finding the generating function $\sum _{\l} ^{}t^{\v\l}$
summed over all partitions $\l$ which
are contained in an $a\times b$ rectangle. As is well-known
(cf\., e.g\., \cite{14, Prop.~1.3.19}), this 
is the $t$-binomial coefficient
$\left[\smallmatrix {a+b}\\ {b} \endsmallmatrix\right]_t$. 
Thus, we obtain the expression (\qtCat).

The claim about $C_n(q,1)$ is the content of Theorem~4.4. For the
proof of the claim about $C_n(1,t)$ we use the correspondence of
Section~2 between ideals $\i$ and partitions $\l$ contained in
$(n,n-1,\dots,1)$. Under this correspondence, the dimension of the
ideal $\i$ is converted into the size $\v\l$ of the partition. If we
rotate the Ferrers diagram of $\l$ by $45^\circ$ in the negative
direction, then the border of the Ferrers diagram forms a Dyck path,
the area of which is exactly equal to $(n+1)^2/2-\v\l$.
\endemo

We are now going to investigate extremal properties, with respect to
dimension and class of nilpotence, of $ad$-nilpotent ideals. First,
we fix the class of nilpotence to $k$, say, 
and ask what the possible dimensions of
ideals with class of nilpotence $k$ is. I.e., the task is to determine
the minimal and maximal possible dimension of an ideal under the assumption  
that its class of nilpotence is $k$. Let us denote the minimal
possible dimension by $\dmin nk$ and the maximal possible dimension by
$\dmax nk$. In terms of our $(q,t)$-Catalan number
$C_n(q,t)$, we ask for the minimal and maximal degree in the variable
$t$ among the terms which have degree $k$ in $q$.

\proclaim{Theorem 6.2}
We have
$$
\dmin nk = {k+1 \choose 2}+(k-1)(n-k)\,,\tag\eqdmin$$
and
$$\dmax nk = \binom {n+1}2 - (n+1)\fl{\frac {n+1} {k+1}}+(k+1)\binom
{\fl{(n+1)/(k+1)}+1} 2.\tag\eqdmax
$$

\endproclaim
\demo{Proof}
In view of the correspondence between ideals and partitions given in
Section~2, determining $\dmin nk$ amounts to finding the partitions
$\l= (\l_1,\l_2,\ldots,\l_n)$ contained in the staircase
$(n,n-1,\dots,1)$ with minimal size $\v\l$ under the
condition that $n(\l)=k$ (i.e., under the condition that the algorithm
of Proposition~3.2 outputs $k$ for $\l$). It is easily seen that under
this assumption we must have $\l_1\ge k$.

Formula (\eqdmin) yields $\dmin n1 = 1$. On the other hand, we have 
$n(1,0,\ldots,0) = 1$ independent of the number of zeroes, since
 the dimension of the ideal associated to  $(1,0,\ldots,0)$ is
always~1. Let us now prove the formula by induction on $k$. We
begin by observing that if
$\l = (\l_1,\l_2,\ldots,\l_n)$ is a partition with $n(\l) = k+1$ and $|\l|$ minimal, then, applying Proposition~3.2, one deduces that 
$|(\l_{n+2-\l_1},\ldots,\l_n)|$ has to be minimal, too. 
Thus, we have to fix $\l_1$ between $k+1$ and $n$, and then take a 
partition $\eta$ which realizes $\dmin{n+2-\l_1}{k}$, 
which we then complete to a partition realizing
$\dmin{n}{k+1}$ in the following way: if $\eta=(s,\ldots)$, then
we take $\l$ as
$(\l_1,s^{n-\l_1},\ldots)$. We thus have
$$
\dmin{n}{k+1} = \min_{k+1 \leq r \leq n} \, \min_{k \leq s \leq r-1}
\{ r + (n-r)s + \dmin{r-1}{k}\}\,. 
$$
By induction, $\dmin{s}{k} = {k+1 \choose 2} + (s-1)(k-1)$, and so
$$
\dmin{n}{k+1} = \min_{k+1 \leq r \leq n} \, \min_{k \leq s \leq r-1}
\left\{ r + (n - r)s + \tfrac{k(k+1)}  2 + (r-1-k)(k-1)\right\}\,. 
$$
Since $n \geq r$ and $k \geq 1$, the minimum on $s$ is reached for $s
= k$, and so 
$$
\dmin{n}{k+1} = \min_{k+1 \leq r \leq n} \left\{ r + nk - rk + 
\tfrac{k(k+1)} 2 + r(k-1) - (k^2-1)\right\}\,.
$$
Thus,
$$
\dmin{n}{k+1} = \min_{k+1 \leq r \leq n} \left\{nk + 
\tfrac{k(k+1)}  2 - k^2 + 1 \right\} = {k+2 \choose 2} + k(n-k-1)\,,
$$
which is what we wanted to prove.

\medskip
Now we consider the case of the maximum.
We want to find partitions $\l$, contained in the staircase
$(n,n-1,\dots,1)$, 
with exactly $k$ outer corners on the antidiagonal  $x+y=n+1$, such
that their size $\v{\l}$ is maximal.

\midinsert
\vskip10pt
\vbox{
$$
\SPfad(0,0),22\endSPfad
\Pfad(0,2),11222111212222\endPfad
\Pfad(0,2),22222222111111\endPfad
\SPfad(6,10),1111\endSPfad
\Label\o{\mu_0}(8,10)
\Label\r{\mu_0}(6,8)
\Label\ro{\mu_1}(5,6)
\Label\lo{\mu_1}(5,5)
\Label\ru{\mu_2}(3,5)
\Label\r{\mu_2}(2,3)
\Label\u{\mu_3}(1,2)
\Label\l{\mu_3}(0,1)
\DuennPunkt(0,0)
\DuennPunkt(0,1)
\DuennPunkt(0,2)
\DuennPunkt(1,2)
\DuennPunkt(2,2)
\DuennPunkt(2,3)
\DuennPunkt(2,4)
\DuennPunkt(2,5)
\DuennPunkt(3,5)
\DuennPunkt(4,5)
\DuennPunkt(5,5)
\DuennPunkt(5,6)
\DuennPunkt(6,6)
\DuennPunkt(6,7)
\DuennPunkt(6,8)
\DuennPunkt(6,9)
\DuennPunkt(6,10)
\DuennPunkt(7,10)
\DuennPunkt(8,10)
\DuennPunkt(9,10)
\DuennPunkt(10,10)
\thinlines
\Diagonale(0,0){10}
\Label\r{x+y=n+1}(11,9)
\DickPunkt(2,2)
\DickPunkt(5,5)
\DickPunkt(6,6)
\hskip7cm
$$
\centerline{\eightpoint Figure 4}
}
\vskip10pt
\endinsert

Let us consider such a partition $\l$, see  
Figure~4 for an example, in which $n=9$ and $k=3$, the partition
being $(6,6,6,6,5,2,2,2,0)$. Let us write
$\mu_0=n+1-\l_1$, $\mu_1=\l_1-\l_2$, \dots,
$\mu_{k-1}=\l_{k-1}-\l_k$, and $\mu_{k}=\l_k$. See Figure~4 for the
geometric meaning of these quantities, where $\mu_0=4$, $\mu_1=1$,
$\mu_2=3$, $\mu_3=2$. Clearly, $n+1=\mu_0+\mu_1+\dots+\mu_k$, i.e.,
$\mu=(\mu_0,\mu_1,\dots,\mu_k)$ is a composition of $n+1$ with exactly
$k+1$ parts. Moreover, the composition $\mu$ determines the partition
$\l$ uniquely. Therefore we may as well encode a partition with all
its outer corners on the antidiagonal $x+y=n+1$ by the corresponding
composition $\mu$. 

We will base our argument on the following two easily verified facts:

\smallskip
{\it Fact 1}. Let $\l$ be a partition (with all outer corners on
$x+y=n+1$) with corresponding composition $\mu$. Let $\mu'$ be a
composition which arose from $\mu$ by permuting the parts. Then the
partition corresponding to $\mu'$ has the same size as $\l$.

{\it Fact 2}. Let $\l$ be a partition with exactly $k$ outer corners,
all of them on $x+y=n+1$, such that its size is maximal with respect
to such partitions. Then, in the corresponding composition
$\mu=(\mu_0,\mu_1,\dots,\mu_k)$, we have $\v{\mu_{i-1}-\mu_{i}}\le1$ for
$i=1,2,\dots,k$. 

\smallskip
Both facts combined say that, in order to find a partition 
with exactly $k$ outer corners,
all of them on $x+y=n+1$, such that its size is maximal, we are
looking for a partition whose corresponding composition $\mu$ has the
property that {\it any two} of its parts differ by at most 1. Thus,
$\mu$ is a composition with at most two different parts, one being
$\fl{(n+1)/(k+1)}$, and the other being $\cl{(n+1)/(k+1)}$. Clearly,
the latter must appear $n+1-(k+1)\fl{(n+1)/(k+1)}$ times, whereas
the former must appear $(k-n+(k+1)\fl{(n+1)/(k+1)})$ times. If one
does the required algebra, then one obtains that the size of such a
partition is exactly the expression on the right hand side of (\eqdmax).
\qed
\enddemo

Now
we fix the dimension to $A$, say,
and ask what the possible classes of nilpotence of
ideals with dimension $A$ is. I.e., now the task is to determine
the minimal and maximal possible classes of nilpotence 
of an ideal under the assumption  
that its dimension is $A$. Let us denote the minimal
possible class of nilpotence 
by $\Dmin nk$ and the maximal possible class of nilpotence by
$\Dmax nk$. In terms of our $(q,t)$-Catalan number
$C_n(q,t)$, we ask for the minimal and maximal degree in the variable
$q$ among the terms which have degree $A$ in $t$.

\proclaim{Theorem 6.3}
We have
$$
\Dmin nA = \min\bigg\{k:\binom {n+1}2 
- (n+1)\fl{\frac {n+1} {k+1}}+(k+1)\binom
{\fl{\tfrac {n+1}{k+1}}+1}2\ge A\bigg\}\,,\tag\eqDmin$$
and
$$\Dmax nA = \fl{n+\tfrac {3} {2}-
\tfrac {1} {2}\sqrt{4n^2+4n+9-8A}}\,.\tag\eqDmax
$$

\endproclaim
\demo{Proof}
Let us first consider the maximum.
We denote the expression on the right hand side of (\eqdmin) by $m(k)$.
Furthermore let $K_0=\max\{k:m(k)\le A\}$. Equivalently, we have
$$m(K_0)\le A< m(K_0+1).\tag\ineqm$$
It is obvious that $\Dmax
nA\le K_0$. We would like to prove equality, because that yields
immediately (\eqDmax) upon a straightforward calculation.

In order to establish $\Dmax nA= K_0$, we consider the partition
$\l_0=(K_0,(K_0-1)^{n-K_0+1},K_0-2,\dots,2,1)$, which realizes the
minimum in $(\eqdmin)$ (with $k=K_0$), i.e.,
$\v{\l_0}=m(K_0)$. It is the lower bound of the interval
$$\multline 
[(K_0,(K_0-1)^{n-K_0+1},K_0-2,\dots,2,1),
(K_0^{n+1-K_0},K_0-1,K_0-2,\dots,2,1)]\\
=[\l^m_{1,2,\dots,K_0},\l^M_{1,2,\dots,K_0}]
\endmultline\tag\eqintm$$
in the decomposition guaranteed by Proposition~4.1. Recall
that all partitions $\l$ in this interval satisfy
$n(\l)=K_0$.  The size of $\l^M_{1,2,\dots,K_0}$ is equal to 
$\binom {K_0}2+(n+K_0-1)K_0$, which is exactly $m(K_0+1)-1$. Hence,
because of (\ineqm), we will be able to find a partition in the
interval (\eqintm) with size $A$. This establishes (\eqDmax).

\medskip
Now we turn to the minimum. Although the idea is analogous, the
details are more elaborate.

We denote the expression on the right hand side of (\eqdmax) by $M(k)$.
Furthermore let $K_1=\min\{k:M(k)\ge A\}$. Equivalently, we have
$$M(K_1)\ge A> M(K_1-1).\tag\ineqM$$
It is obvious that $\Dmin
nA\ge K_1$. We would like to prove equality.

In order to establish $\Dmin nA= K_1$, we consider the partition
$$\multline
\l_1=\Big(\(n+1-\cl{\tfrac {n+1}{K_1+1}}\)^{\cl{(n+1)/(K_1+1)}},
\(n+1-2\cl{\tfrac
{n+1}{K_1+1}}\)^{\cl{(n+1)/(K_1+1)}},\\
\dots,
\(2\fl{\tfrac {n+1}{K_1+1}}\)^{\fl{(n+1)/(K_1+1)}},
\(\fl{\tfrac {n+1}{K_1+1}}\)^{\fl{(n+1)/(K_1+1)}}\Big),
\endmultline\tag\eqstufen$$ 
(to be precise, the partition
corresponding to the composition $\big(\clkl{\tfrac {n+1}{K_1+1}}^a,
\flkl{\tfrac {n+1}{K_1+1}}^b\big)$, where we have abbreviated
$a=n+1-(K_1+1)\flkl{\tfrac {n+1}{K_1+1}}$ 
and $b=K_1-n+(K_1+1)\flkl{\tfrac {n+1}{K_1+1}}$; see the second part
of the proof of Theorem~6.2),
which realizes the
maximum in (\eqdmax) (with $k=K_1$), i.e.,
$\v{\l_1}=M(K_1)$. It is the upper bound of the interval
$$[\l^m_{(i_1,\dots,i_{K_1})},
\l^M_{(i_1,\dots,i_{K_1})}],
\tag\eqintM$$
where
$$(i_1,\dots,i_{K_1})=
\Big(\fl{\tfrac {n+1}{K_1+1}},2\fl{\tfrac {n+1}{K_1+1}},\dots,
n+1-2\cl{\tfrac
{n+1}{K_1+1}},n+1-\cl{\tfrac {n+1}{K_1+1}}\Big),$$
in the decomposition guaranteed by Proposition~4.1. Recall
that all partitions $\l$ in this interval satisfy
$n(\l)=K_1$.  As a moderately tedious computation shows,
the size of $\l^m_{(i_1,\dots,i_{K_1})}$ is equal to 
$$\cases \binom {n+1}2 - (3n+4)\fl{\frac {n+1} {K_1+1}}+(3K_1+5)\binom
{\fl{(n+1)/(K_1+1)}+1} 2&\text{if }(K_1+1)\nmid(n+1),\\
\binom {n+1}2 - (3n+5){\frac {n+1} {K_1+1}}+(3K_1+5)\binom
{{(n+1)/(K_1+1)}+1} 2&\text{if }(K_1+1)\mid(n+1).
\endcases\tag\eqintmin$$
If we are able to establish that this value is less than or equal to
$M(K_1-1)+1$, then, because of (\ineqM), we will be able to find a
partition in the 
interval (\eqintM) with size $A$. This would establish (\eqDmin). 

In fact, it turns out that the preceding claim is only true for
$K_1>1$. Hence, we will treat the case of $K_1=1$ separately at the
end of the proof.

Let $K_1>1$. First, the claim is easily verified directly for $K_1=n$
(in which case the expression in the second line of (\eqintmin) has to
be used).
Second, we verify our claim for $n=1,2,\dots,6$. This is
readily done with the help of a computer. (It can even be done by
hand.) Thus, in
the sequel, we may assume that $K_1\le n-1$ and $n\ge 7$.

Since the expression in the second line in (\eqintmin) is
smaller than the expression in the first line, it suffices to
prove that the expression in first line is less or equal to
$M(K_1-1)+1$. That is, we must show
$$\multline \binom {n+1}2 - (3n+4)\fl{\frac {n+1} {K_1+1}}+(3K_1+5)\binom
{\fl{(n+1)/(K_1+1)}+1} 2\\
\le
\binom {n+1}2 - (n+1)\fl{\frac {n+1} {K_1}}+K_1\binom
{\fl{(n+1)/K_1}+1} 2 +1
\endmultline$$
This inequality is equivalent to
$$\multline 0\le (3n+4)\fl{\frac {n+1} {K_1+1}}
-(3K_1+5)\binom
{\fl{(n+1)/(K_1+1)}+1} 2\\
+K_1\binom
{\fl{(n+1)/K_1}+1} 2 - (n+1)\fl{\frac {n+1} {K_1}}
+1.
\endmultline\tag\ineqA$$
It should be observed that, as long as $K_1$ is between $n/2+1$ and
$n-1$, the
right hand side of (\ineqA) is linear and monotone decreasing in
$K_1$. On the other hand, it is trivially true for $K_1=n-1$. Hence,
it is true for all $K_1\ge n/2+1$. This allows us to assume $K_1\le
(n+1)/2$ from now on.

The expression on the right hand side of (\ineqA) is quadratic in
$\fl{(n+1)/K_1}$, with the minimum of the quadratic polynomial at
$(n+1)/K_1-1/2$. Thus, if we would be able to prove (\ineqA) with
$\fl{(n+1)/K_1}$ replaced by $(n+1)/K_1-1/2$, the original inequality
would be established. Similarly, the right hand side of (\ineqA) is
quadratic in  $\fl{(n+1)/(K_1+1)}$. Therefore, depending on whether 
$\fl{(n+1)/(K_1+1)}$ is to the right or to the left of the maximum of the
corresponding quadratic polynomial, it suffices to prove (\ineqA) with 
$\fl{(n+1)/(K_1+1)}$ replaced by $(n+1)/(K_1+1)$, respectively
$(n-K_1+1)/(K_1+1)$. 

In summary, we will be done once we have established the inequalities
$$0\le
 \frac {\eqalign{
&(-8 n-12 {K_1}-9 {K_1}^2-4 n^2-4-4 {K_1} n^2\cr
&\hskip4cm
-16 {K_1} n+8 {K_1}^2 n^2-8 {K_1}^3 n-10 {K_1}^3-{K_1}^4 )}
}{8{K_1}({K_1}+1)^2},\tag\ineqa$$
corresponding to replacing $\fl{(n+1)/K_1}$ by $(n+1)/K_1-1/2$ and
$\fl{(n+1)/(K_1+1)}$ by $(n+1)/(K_1+1)$ in (\ineqA), and
$$0\le 
 \frac {\eqalign{
&(-8 n-12 {K_1}+19 {K_1}^2-4 n^2-4-4 {K_1} n^2-16 {K_1} n\cr
&\hskip4cm+8 {K_1}^2 n^2+16 {K_1}^2 n-8 {K_1}^3 n-6 {K_1}^3-{K_1}^4)}} 
{8{K_1}({K_1}+1)^2},\tag\ineqb$$
corresponding to replacing $\fl{(n+1)/K_1}$ by $(n+1)/K_1-1/2$ and
$\fl{(n+1)/(K_1+1)}$ by $(n-K_1+1)/(K_1+1)$ in (\ineqA).

We concentrate on the proof of (\ineqa). The proof of (\ineqb) is
similar. 

First of all, it can be verified directly that (\ineqa) is true for
$K_1=2$ and $n\ge7$. Therefore, from now on, we may assume $K_1\ge3$.
Next we differentiate the expression on the right hand side of
(\ineqa) with respect to $K_1$, thus obtaining
$$ -\frac {\eqalign{
&(n^2(8 {K_1}^3 -16 {K_1}^2-12 {K_1}-4)+n(16 {K_1}^3-32 {K_1}^2-24 {K_1}-8)\cr
&\hskip5cm+
{K_1}^5+3 {K_1}^4+11 {K_1}^3-15 {K_1}^2-12 {K_1})}} {8{K_1}^2({K_1}+1)^3}.$$
This is most evidently negative for $K_1\ge3$. Hence, the right hand
side of (\ineqa) is monotone decreasing for $K_1\ge3$. Thus, if we are
able to verify (\ineqa) for the maximal $K_1$ that we are considering,
i.e., $K_1=(n+1)/2$, then (\ineqa) is established for all $K_1$
between $3$ and $(n+1)/2$. Now, if we substitute $K_1=(n+1)/2$ into
(\ineqa), we obtain
$$\frac {(n+1)(15n^2-70n-217)} {16(n+3)^2},$$
which is positive for $n\ge7$, as desired.

Finally, we treat the case $K_1=1$. In that case, because of (\ineqM),
our given size $A$ must satisfy $1\le A\le \fl{(n+1)/2}\cl{(n+1)/2}$.
If $A\ge n$, then there is a partition in the interval
$[I^m_{\fl{(n+1)/2}},I^M_{\fl{(n+1)/2}}]$ with size $A$. Otherwise,
there is a partition in the interval $[I^m_{1},I^M_{1}]$ with size $A$.

This completes the proof of the theorem. \qed
\enddemo

%\par\newpage
\enddocument
\bye